\documentclass{agtart_a}
\pdfoutput=1
\usepackage{graphicx}

\title{The $\mathrm{FA}_n$ Conjecture for Coxeter groups}

\author{Angela Kubena Barnhill}
\givenname{Angela Kubena}
\surname{Barnhill}
\address{Department of Mathematics\\The Ohio State University\\\newline
231 West 18th Avenue\\Columbus, Ohio  43210\\USA}
\email{abarnhill@math.ohio-state.edu}
\urladdr{}

\volumenumber{6}
\issuenumber{}
\publicationyear{2006}
\papernumber{73}
\startpage{2117}
\endpage{2150}

\doi{}
\MR{}
\Zbl{}


\received{25 September 2005}
\revised{}
\accepted{6 March 2006}
\published{19 November 2006}
\publishedonline{19 November 2006}
\proposed{}
\seconded{}
\corresponding{}
\editor{CPR}
\version{}

\arxivreference{math.GR/0509439}

\AtBeginDocument{\let\bar\wbar\let\tilde\wtilde\let\hat\what}

\makeatletter
\def\cnewtheorem#1[#2]#3{\newtheorem{#1}{#3}[section]
\expandafter\let\csname c@#1\endcsname\c@theorem}
\makeatother

\theoremstyle{plain}

\newtheorem{theorem}{Theorem}[section]
\cnewtheorem{prop}[theorem]{Proposition}
\cnewtheorem{cor}[theorem]{Corollary}
\cnewtheorem{conj}[theorem]{Conjecture}

\newtheorem{CATOconj}{CAT(0) Conjecture}
\makeautorefname{CATOconj}{CAT(0) Conjecture}

\newtheorem*{FAnconj}{Maximal $\FA_n$ Conjecture}
\cnewtheorem{lemma}[theorem]{Lemma}

\theoremstyle{definition}
\cnewtheorem{defn}[theorem]{Definition}
\makeautorefname{defn}{Definition}
\newtheorem*{notation}{Notation}
\cnewtheorem{eg}[theorem]{Example}
\cnewtheorem{remark}[theorem]{Remark}
\cnewtheorem{note}[theorem]{Note}
\newtheorem*{acknow}{Acknowledgements}

\def\RR{\mathbb R}
\def\NN{\mathbb N}
\def\ZZ{\mathbb Z}
\def\QQ{\mathbb Q}
\def\HH{\mathbb H}
\def\EE{\mathbb E}
\def\SS{\mathbb S}
\def\CC{\mathbb C}

\def\implies{\Rightarrow}
\def\cS{\Sigma}
\def\cN{\mathcal N}
\def\cT{\mathcal T}
\def\st{:}

\DeclareMathOperator{\FA}{FA}

\DeclareMathOperator{\CAT}{CAT}
\DeclareMathOperator{\SL}{SL}
\DeclareMathOperator{\GL}{GL}
\DeclareMathOperator{\Isom}{Isom}

\DeclareMathOperator{\Fix}{Fix}

\DeclareMathOperator{\alt}{alt}

\begin{document}

\begin{asciiabstract}
We study global fixed points for actions of Coxeter groups on
nonpositively curved singular spaces.  In particular, we consider
property FA_n, an analogue of Serre's property FA for actions on
CAT(0) complexes.  Property FA_n has implications for irreducible
representations and complex of groups decompositions.  In this paper,
we give a specific condition on Coxeter presentations that implies FA
_n and show that this condition is in fact equivalent to FA _n for n=1	
and 2.  As part of the proof, we compute the Gersten-Stallings angles
between special subgroups of Coxeter groups.
\end{asciiabstract}

\begin{htmlabstract}
We study global fixed points for actions of Coxeter groups on
nonpositively curved singular spaces.  In particular, we consider
property FA<sub>n</sub>, an analogue of Serre's property FA for
actions on CAT(0) complexes.  Property FA<sub>n</sub> has
implications for irreducible representations and complex of groups
decompositions.  In this paper, we give a specific condition on
Coxeter presentations that implies FA<sub>n</sub> and show that this
condition is in fact equivalent to FA<sub>n</sub> for n=1 and 2.
As part of the proof, we compute the Gersten&ndash;Stallings angles between
special subgroups of Coxeter groups.
\end{htmlabstract}

\begin{abstract}
We study global fixed points for actions of Coxeter groups on
nonpositively curved singular spaces.  In particular, we consider
property $\mathrm{FA}_n$, an analogue of Serre's property FA for
actions on $\mathrm{CAT}(0)$ complexes.  Property $\mathrm{FA}_n$ has
implications for irreducible representations and complex of groups
decompositions.  In this paper, we give a specific condition on
Coxeter presentations that implies $\mathrm{FA}_n$ and show that this
condition is in fact equivalent to $\mathrm{FA}_n$ for $n=1$ and 2.
As part of the proof, we compute the Gersten--Stallings angles between
special subgroups of Coxeter groups.
\end{abstract}

\maketitle

\section{Introduction}

A {\em Coxeter group} is a group $W$ that has a presentation of the form \[W=\left< S\mid(s_is_j)^{m_{ij}} = 1\right>\] where $ m_{ij}=m_{ji} \in \NN \cup \{\infty \}$ and $m_{ij} = 1$ if and only if $i=j$.  Recall that a {\em $\CAT (0)$ space} is a complete geodesic space which is nonpositively curved in the metric sense, ie, its geodesic triangles are no fatter than their Euclidean counterparts (see \fullref{CAT0basics}).  We will consider isometric actions of Coxeter groups on $\CAT (0)$ spaces.

A fundamental notion of Bass--Serre theory is Serre's {\em property $\FA$}.  A group $G$ has property FA if every $G$--action on every simplicial tree is {\em trivial}, ie, has a global fixed point.  Such groups are ``rigid" in the following sense: they do not split nontrivially as amalgamated free products or HNN extensions and all their irreducible GL$_2( \CC )$--representations have algebraic integer traces.  In addition to finite groups, Serre proved in \cite{SerreTrees} that $\SL _3 (\ZZ)$ and Coxeter groups with every $m_{ij}<\infty$ have FA.

A generalization of Serre's property is {\em property $\FA _n$}, as defined by Farb in \cite{FarbHelly}.  A group has $\FA _n$ if every action of it by cellular isometries on every $\CAT (0)$ \hbox{$n$--complex} is trivial.  By $\CAT (0)$ $n$--complex, we mean a $\CAT (0)$ cell complex of piecewise constant curvature with only finitely many isometry types of cells (see \fullref{CAT0basics}).  We emphasize that the actions are not assumed to be cocompact, properly discontinuous or faithful, and that the spaces are not assumed to be locally finite.  Note that FA$_1$ is equivalent to FA.

Properties $\FA _n$ and FA$_m$ are distinct for $n \neq m$.  For example, Farb proved in \cite{FarbHelly} that $\SL_n(\ZZ[{1}/{p}])$ has FA$_{n-2}$.  However, $\SL_n(\ZZ[1 / p ])$ does not have FA$_{n-1}$ since it acts nontrivially on the affine building for $\SL_n(\QQ _p)$.

As with property FA (see \cite{SerreTrees}), groups with $\FA _n$ have certain strong properties.
\begin{enumerate}
\item  If $\Gamma$ has property $\FA_n$ then $\Gamma$ does not split nontrivially as a nonpositively curved {\em $n$--complex of groups} in the sense of Gersten--Stallings, Haefliger and Corson (see \cite{StallingsTriangles,Haefliger,corson}).
\item  Suppose $\Gamma$ has property $\FA_n$.  Let $\rho \co \Gamma \rightarrow \GL _{n+1}(K)$ be any representation of degree $n+1$ over a field $K$.  Then, the eigenvalues of each of the matrices in $\rho(\Gamma)$ are integral.  In particular, they are algebraic integers if char$(K)=0$ and are roots of unity if char$(K) > 0$.  (As in the tree case, this follows from studying induced actions on the Bruhat--Tits buildings for SL$_{n+1}(\QQ_p)$, which are $\CAT (0)$ for all primes $p$.)  In the language of Bass \cite{Bass80}, $\Gamma$ is thus of {\em integral $(n+1)$--representation type}.  Consequently, there are only finitely many conjugacy classes of irreducible representations of $\Gamma$ into $\GL _{n+1}$(K) for any algebraically closed field K (see Farb \cite{FarbHelly}).
\end{enumerate}

The following results are known about property $\FA _n$ for a Coxeter group $W$:
\begin{enumerate}
\item (Serre \cite{SerreTrees}) If every $m_{ij}$ is finite, then $W$ has property FA.
\item (Farb \cite{FarbHelly}) If $W$ is a discrete group generated by reflections in the sides of a compact Euclidean or hyperbolic $n$--simplex, then $W$ has property FA$_{n-1}$ but  does not have property $\FA _n$.
\end{enumerate}
In this paper, we generalize these results.  We first consider natural conditions on Coxeter groups that imply property $\FA _n$.

For $T \subset S$, let $W_T$ denote the subgroup of $W$ generated by $T$.  It is well-known (see, for example, Bourbaki \cite{bourbaki}) that $(W_T,T)$ is a Coxeter system with a Coxeter presentation that is induced from the presentation for $W$.  The group $W_T$ is a {\em special subgroup} with {\em rank} equal to the size of $T$.

Applying techniques from \cite{FarbHelly}, we prove the following theorem by considering the combinatorics of fixed sets of finite special subgroups.
\begin{theorem}
\label{thmhelly}
Let $(W,S)$ be a Coxeter system.  If every special subgroup of $W$ of rank at most $n+1$ is finite, then W has property $\FA _n$.
\end{theorem}

Suppose a group acts nontrivially on an $n$--dimensional $\CAT (0)$ space $X$.  Then for $m\geq n$ it acts nontrivially on the $m$--dimensional $\CAT (0)$ space $X \times \RR^{m-n}$.  So, for $m\geq n$ we have FA$_m \Rightarrow $ $\FA _n$.  In other words, for every group $G$ that acts nontrivially on some finite-dimensional $\CAT (0)$ complex, there an integer $n$ such that $G$ has $\FA _m$ if and only if $m<n$.  This $n$ is the smallest dimension of a $\CAT (0)$ complex on which $G$ acts nontrivially.  As formulated in the following conjecture, we suspect that \fullref{thmhelly} gives this bound for Coxeter groups.
\begin{conj}
\label{mainconj}
 Let $(W,S)$ be a Coxeter system.  The following are equivalent:
\begin{enumerate}
\item[(i)]  The group $W$ has property $\FA _n$.
\item[(ii)]  Every special subgroup of $W$ of rank at most $n+1$ is finite.
\item[(iii)] For all $0<m\leq n$, the group $W$ does not split nontrivially as a nonpositively curved $m$--simplex of special subgroups.
\end{enumerate}
\end{conj}

As noted by Mihalik and Tschantz in \cite{MihalikTschantz}, \fullref{mainconj} is known for $n=1$.  In this paper, we reduce the proof of \fullref{mainconj} in general to proving that spaces arising from certain simplex of groups decompositions of $W$ are $\CAT (0)$.  The $\CAT (0)$ Conjecture (\fullref{construct}) posits that these spaces are indeed $\CAT (0)$.  In \fullref{Dim2Section}, we prove the $\CAT (0)$ Conjecture in dimension 2 by computing the Gersten--Stallings angles between special subgroups of Coxeter groups.  This implies the following theorem.

\begin{theorem}
\label{thm12}
\fullref{mainconj} holds for $n\leq 2$.
\end{theorem}

P Caprace observed that \fullref{mainconj} holds for all $n$ if and only if it holds for $n\leq 8$ (see \fullref{caprace}).  To prove \fullref{mainconj} in general, it therefore remains only to show that it holds for $3\leq n\leq 8$.

In \fullref{section:MaxFAn}, we study the maximal $\FA_n$ subgroups of Coxeter groups.  Special subgroups satisfying the condition of \fullref{thmhelly} are natural candidates, and we posit the following:

\begin{FAnconj}
A subgroup $H \subset W$ is maximal $\FA _n$ if and only if $H = wAw^{-1}$ for some maximal $\FA_n$ special subgroup $A$ of $W$ and $w \in W$.
\end{FAnconj}

This has been shown for $n=1$ by Mihalik and Tschantz in \cite{MihalikTschantz}.  In \fullref{section:MaxFAn}, we prove the following:

\begin{theorem} \label{maxFAnfromCAT0}
$\CAT(0)$ Conjecture $\implies$ Maximal $\FA_n$ Conjecture.
\end{theorem}

In particular, since the $\CAT (0)$ Conjecture holds in dimension 2, the Maximal $\FA_2$ Conjecture holds as well.

In Sections $2$--$4$, we recall important results about Coxeter groups, $\CAT (0)$ spaces, and complexes of groups.  We discuss in Section 5 background and techniques related to property $\FA _n$.  We present and prove our main results in Sections $6$--$8$.

Recently, Luis Paris has independently discovered results similar to those in Section 6.

\begin{acknow}
I would first like to thank my thesis advisor Benson Farb for his mathematical guidance, time and patience.  His ideas, including asking the questions that led me to this project, and his help in the preparation of this paper have been invaluable.  I am also grateful to Mike Davis for his time and insights.  Many thanks to Roger Alperin, Ruth Charney, Chris Hruska and Kevin Whyte for their helpful conversations, to Mark Behrens, Daniel Biss, Jon McCammond, Mike Mihalik and Luis Paris for their suggestions, and to Alissa Crans, Pallavi Dani, Ian Leary and Anne Thomas for their helpful comments on this paper.  I would also like to thank Pierre--Emmanuel Caprace for his comments on the $\FA_n$ Conjecture for large values of $n$.  Finally, I thank the referee for many useful suggestions.
\end{acknow}

\section{Coxeter groups}
We briefly recall key definitions and results on Coxeter groups that we will need.  See Bourbaki \cite{bourbaki}, Davis \cite{Davisbook}  or Humphreys \cite{Humphreys} for further details.

\subsection{Definitions} \label{Coxeterdefs}
Let $S$ be a finite set.  A {\em Coxeter matrix} on $S$ is a symmetric
 $|S| \times |S|$ matrix $M$ with entries in $\NN \cup \{\infty\}$
 such that each diagonal entry is 1 and each off-diagonal entry is at
 least 2.  Associated to $M$ is a group $W$ with presentation
 $W=\left\langle S \mid (st)^{m_{st}} = 1 \right\rangle$, where the
 relation $(st)^{m_{st}}=1$ is omitted if $m_{st}=\infty$.  The pair
 $(W,S)$ is a {\em Coxeter system}, $W$ is a {\em Coxeter group}, and
 the group presentation is a {\em Coxeter presentation} for $W$.  The
 $W$--conjugates of elements of $S$ are called {\em reflections}.

Given a Coxeter system $(W,S)$, the {\em Coxeter diagram} $\Gamma$ associated to $(W,S)$ is a labeled graph with vertex set $S$ and with an edge labeled $m_{st}$ connecting $s$ to $ t$ if and only if $m_{st} \neq 1,2$.  Note that the Coxeter diagram encodes the same information as the Coxeter matrix.  A Coxeter system $(W,S)$ is {\em irreducible} if its Coxeter diagram is a connected graph.

\subsection{Reduced words and the word problem}
Given a Coxeter system $(W,S)$ and $w\in W$, we denote by ${\ell} (w)$ the {\em length} of $w$, which is given by \[ \ell (w)= \min\{k \st w=s_{1}s_{2}\ldots s_{k} \textrm{ for some }s_{1},s_{2},\ldots,s_{k} \in S\}.\]  An expression for $w$ which achieves its length $\ell (w)$ is called a {\em reduced} or {\em geodesic} expression for $w$.  Note that this notion depends on the choice of Coxeter generating set $S$.

We denote by $\widehat{s_i}$ the omission of $s_i$ from an expression.  Below is a standard characterization of Coxeter groups (see, for example, \cite{BrownBuildings,Davisbook,Humphreys}).

\begin{theorem}[Deletion and Exchange Conditions] \label{deletion} Let $W$ be a group generated by a set $S$ of involutions.  The following are equivalent:
\begin{enumerate}
\item[(i)] $(W,S)$ is a Coxeter system.
\item[(ii)] {\em(Deletion Condition)} For all $w \in W$, if $\ell (w) < k$ and $w=s_{1}s_{2}\ldots s_{k}$ for some generators $s_1, s_2, \ldots s_k \in S$, then there exist indices $1\leq j<l\leq k$ such that $w=s_{1}s_{2} \ldots \widehat{ s_{j}} \ldots {\widehat s_{l}} \ldots s_{k}$.
\item[(iii)] {\em(Strong Exchange Condition)} Let $w \in W$ and let $w=s_1s_2\ldots s_k$ ($s_i \in S$) be an expression for $w$.  If a reflection $r$ in $W$ satisfies \hbox{$\ell (rw) < \ell(w)$}, then there is an index $i$ for which \hbox{$w = rs_1\ldots \widehat{s_i} \ldots s_k$}.  Moreover, if $k=\ell (w)$, then $i$ is unique.
\end{enumerate}
\end{theorem}

An immediate consequence of the Deletion Condition is the following:
\begin{cor}\label{geodesicwalls}
Let $(W,S)$ be a Coxeter system and $w \in W$.  Then every unreduced expression for $w$ can be reduced to a geodesic expression for $w$ by omitting an even number of generators.
\end{cor}
We now define a standard set of operations for reducing words in Coxeter groups.
\begin{defn}\label{elemdefn}
Let $M$ be the Coxeter matrix associated to a Coxeter system $(W,S)$.  An {\em elementary $M$--operation} on a word $w$ in the alphabet $S$ is an operation of {\em Type (I)} or {\em Type (II)}, which are defined as:
\begin{itemize}\leftskip 25pt
\item[\bf Type (I)]Delete a subword of the form $ss$ for some $s \in S$.
\item[\bf Type (II)]Replace an alternating subword $sts\ldots$ of length $m_{st}$ for some $s,t\in S$  with the alternating word $tst\ldots$ of length $m_{st}$.
\end{itemize}
Note that elementary $M$--operations do not change the image of the word in $W$.  An {\em $M$--reduction} of a word is a sequence of elementary $M$--operations.  A word is {\em $M$--reduced} if its length cannot be reduced via elementary $M$--operations.
\end{defn}

The {\em word problem} is a fundamental problem in combinatorial and geometric group theory.  The Deletion Condition implies that Coxeter groups have solvable word problem.  In particular, Tits (see \cite{TitsWP,BrownBuildings}) proved the following theorem.

\begin{theorem} \label{titswp}
Let $(W,S)$ be a Coxeter system with associated Coxeter matrix $M$.  Then an expression for $w\in W$ is reduced if and only if it is $M$--reduced.  Moreover, given two reduced expressions for $w$, one can be transformed to the other via a sequence of Type (II) elementary $M$--operations.
\end{theorem}

Since elementary $M$--operations do not increase word length, \fullref{titswp} solves the word problem for Coxeter groups.

\begin{remark} \label{uniquealtprod}
An important special case of \fullref{titswp} is for alternating products.  In particular, for $a,b,\in S$, an alternating product of $a$ and $b$ of length strictly less than $m_{ab}$ is the {\em unique} reduced representative of the corresponding element of $W$.
\end{remark}

\subsection{Parabolic subgroups}
Given a Coxeter group $W$, we will be particularly interested in certain natural subgroups of $W$.

Note that for elements or subsets $A_1, A_2 , \ldots , A_k$ of a group $G$ we denote by\break $\langle A_1,A_2 , \ldots ,A_k\rangle$ the subgroup of $G$ generated by the union of the $A_i$'s.
\begin{defn}
A {\em special} subgroup of $W$ is a subgroup $W_T$ of $W$ given by $W_T = \langle T \rangle$ for some $T \subset S$.  We say $W_\emptyset = 1$ and we define the {\em rank} of $W_T$ to be $|T|$.  The $W$--conjugates of special subgroups of $W$ are called {\em parabolic} subgroups of $W$.
\end{defn}

\begin{theorem}[See, for example, \cite{Humphreys}] \label{specialchar}  Let $(W,S)$ be a Coxeter system with Coxeter matrix $M$, and let $T\subset S$.
\begin{enumerate}
\item[(i)] $(W_T,T)$ is a Coxeter system with Coxeter matrix the submatrix of $M$ corresponding to $T$.
\item[(ii)] If $w \in W_T$, and $w=s_1s_2\ldots s_k$ is a reduced expression for $w$ with respect to $S$, then $s_i \in T$ for all $i$.  In particular, $W_T \cap S=T$.  Moreover, the length of $w$ with respect to $S$ equals the length of $w$ with respect to $T$.
\item[(iii)] $S$ is a minimal generating set for $W$, and more generally, $T$ is a minimal generating set for $W_T$.
\end{enumerate}
\end{theorem}
\eject

Given $T\subset S$, consider the Coxeter system $(W_T,T)$.  Denote by $\Gamma_T$ the associated Coxeter diagram.  Note that by \fullref{specialchar}, this is the labeled subgraph of $\Gamma$ spanned by the vertices corresponding to the elements of $T$.

\begin{remark} \label{specialinter} Given a special subgroup $G\subset W$, we denote by $S_G$  the generators of $G$ as a special subgroup.  By \fullref{specialchar}, we have \hbox{$S_G = G \cap S$} and \hbox{$G=W_{S_G}$}.  In this notation, statement (ii) of \fullref{specialchar} implies the well-known result that for special subgroups $A$ and $B$, we have \hbox{$A\cap B=\langle S_A \cap S_B \rangle$.}
\end{remark}
For $T,U\subset S$, we say that $w$ is {\em(T,U)--reduced} if it is of minimal length in the double coset $W_TwW_U$.  The following is standard (see, for example,  \cite{bourbaki}).
\begin{prop} \label{XX'reduced}
Let $T$ and $U$ be (possibly empty) subsets of $S$.
\begin{enumerate}
\item[(i)] Let $w\in W$. There is a unique $(T,U)$--reduced element $d\in W_TwW_U$.  Moreover, every element $w' \in W_TwW_U$ can be written as \hbox{$w'=xdy $} for some $x \in W_T$ and $y \in W_U$ so that $\ell (w') = \ell (x) + \ell(d) + \ell (y)$.
\item[(ii)] An element $w\in W$ is $(T,U)$--reduced if and only if $\ell (tw) > \ell (w) \ \forall \ t \in T$ and $\ell(wu) > \ell (w) \ \forall \ u \in U$
\end{enumerate}
\end{prop}

Given two parabolic subgroups of $W$, we will be interested in their intersection.  In fact, the intersection of two parabolic subgroups is again a parabolic subgroup.  In particular, we have the following result of Kilmoyer (see \cite{MihalikTschantz}).

\begin{prop}\label{paris}
Suppose $A$ and $B$ are special subgroups of $W$, with corresponding generating sets $S_A, S_B \subset S$.  For $w\in W$, let $d$ be the unique $(S_A,S_B)$--reduced element of $AwB$ and let $a\in A$, $b\in B$ so that $w=adb$.  Then $$A \cap w B w^{-1} =  a\langle S_A \cap dS_B d^{-1} \rangle a^{-1}.$$
\end{prop}

\begin{cor}\label{inclgen}
Suppose $S_A, S_B \subset S$ and let $w \in W$. Let $d$ be the unique
$(S_A,S_B)$--reduced element in $AwB$.  If $A\subset wBw^{-1}$, then
$S_A \subset dS_Bd^{-1}$.  \end{cor}

\begin{proof}[Proof of \fullref{inclgen}]
We have $A = A \cap w B w^{-1}$, so for $a \in A$ and $b\in B$ such that $w=adb$,  we find $A= a \langle S_A \cap dS_B d^{-1} \rangle a^{-1}$ by \fullref{paris}.  Since $a\in A$, we know $a^{-1}Aa=A$ and hence $A=\langle S_A \cap dS_B d^{-1}\rangle$.   Now $S_A$ is a minimal generating set of $A$ by statement (iii) of \fullref{specialchar}, so $A$ cannot be generated by a proper subset of $S_A$.  We therefore find that $S_A = S_A \cap  dS_B d^{-1}$ so in particular $S_A \subset dS_Bd^{-1}$.
\end{proof}

\subsection{Classification}
Let $W$ be a Coxeter group arising as a discrete group generated by reflections in the sides of a compact Euclidean or hyperbolic simplex.  Note that every proper special subgroup of $W$ is finite since every such subgroup stabilizes a point under the discrete and proper action of $W$ on Euclidean or hyperbolic space.  In fact, this property distinguishes these groups (see, for example, \cite{BrownBuildings,Davisbook,Humphreys}).
\begin{theorem} \label{classification}
Let $(W,S)$ be a Coxeter system of rank $n+1$.  Then, every proper special subgroup of $W$ is finite if and only if one of the following holds:
\begin{enumerate}
\item[(i)]  $W$ is finite.
\item[(ii)] $W$ is an irreducible Euclidean or hyperbolic reflection group with fundamental domain a compact simplex.
\end{enumerate}
\end{theorem}

In case (ii), $W$ acts cocompactly on a Euclidean or hyperbolic space of dimension $n$, and we can recover the dihedral angles of the fundamental domain from the Coxeter presentation.  Moreover, the hyperbolic case only occurs for $n\leq 4$.
\begin{remark}\label{lookupfinite}
This is part of the general classification of Coxeter groups.  Irreducible finite Coxeter groups and those described in case (ii) of \fullref{classification} are listed in standard books on Coxeter groups, such as \cite{bourbaki} and \cite{Humphreys}.  In particular, determining whether Coxeter groups (or their special subgroups) are finite reduces to verifying whether their irreducible components appear on the list of irreducible finite Coxeter groups.
\end{remark}

\section{CAT(0) spaces}
We recall a few key facts about $\CAT (0)$ spaces.  See Bridson and Haefliger \cite{BridsonHaefliger} for full details.

\subsection{Definitions, properties and constructions} \label{CAT0basics}

\begin{figure}[ht!] \label{comptri}
\begin{center}
\includegraphics{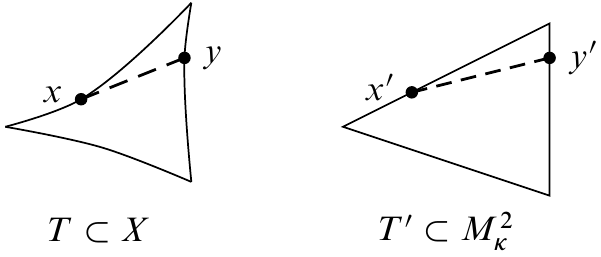}
\caption{Comparison triangle}
\end{center}
\end{figure}

Given $\kappa \in \RR$, let $M^n_{\kappa}$ denote $\EE^n$ (Euclidean $n$--space), $\HH^n_{\kappa}$ (hyperbolic $n$--space of constant curvature $\kappa$) or $\SS^n_{\kappa}$ (the $n$--sphere of constant curvature $\kappa$) as $\kappa$ is 0, negative or positive, respectively.  Denote by $d_{\kappa}$ the distance function on $M^2_{\kappa}$.
\eject
\begin{defn}
Assume that $(X,d)$ is a geodesic metric space.  Let $T \subset X$ be a geodesic triangle.  A {\em $\kappa$--comparison triangle} for $T$ is a triangle $T' \subset M^2_{\kappa}$ with the same edge lengths as those in $T$.  We say $T$ satisfies the {\em CAT($\kappa$) inequality} if for all points $x , y \in T$ with corresponding points $x', y' \in T'$ (see \fullref{comptri}), we have $d(x,y) \leq d_{\kappa} (x',y')$.  The space $X$ is a {\em CAT($\kappa$) space} if every triangle $T\subset X$ (of perimeter at most $\frac {2 \pi}{\sqrt{\kappa}}$ if $\kappa > 0$) satisfies the $\CAT(\kappa)$ inequality.  Note that $\CAT(\kappa)$ implies $\CAT(\kappa ')$ for all $\kappa ' \geq \kappa$.
\end{defn}

\begin{remark}  \label{CAT0props}
We will be most interested in the case $\kappa = 0$.  $\CAT(0)$ spaces are a generalization of complete, simply connected, nonpositively curved Riemannian manifolds.  $\CAT(0)$ spaces can be singular and even locally infinite, but the $\CAT(0)$ condition implies many strong properties, including convexity of the distance function, unique geodesics and contractibility.
\end{remark}

We now briefly discuss a natural construction that yields many interesting examples of $\CAT (0)$ spaces.

An {\em $M_{\kappa}$--polyhedral complex} is a cell complex formed by taking the disjoint union of convex polyhedral cells in $M^n_{\kappa}$ which are then glued along isometric faces.  An $M_{\kappa}$--polyhedral complex is called piecewise-Euclidean, piecewise-hyperbolic or piecewise-spherical if $\kappa$ is 0, -1 or 1, respectively.  (See \cite{BridsonHaefliger} for details.)

Bridson showed in \cite{Bridson91} that if an $M_{\kappa}$--polyhedral complex is composed of only finitely many isometry types of cells, then it is a complete geodesic space with respect to the naturally defined metric.  Note that this condition is often satisfied in cases that arise naturally.  For example, if a group acts cocompactly by isometries on a metric polyhedral complex $X$, then $X$ has only finitely many isometry types of cells.

\begin{defn} \label{linkconddefn}
An $M_{\kappa}$--polyhedral complex $X$ satisfies the {\em link condition} if for every vertex $v \in X$, the link of $v$ (with its natural piecewise-spherical structure) is a $\CAT(1)$ space.
\end{defn}

The importance of this condition is given by the following theorem, due to Gromov, Ballman and Bridson (see \cite{BridsonHaefliger}).

\begin{theorem} \label{linkcond}
Let $X$ be an $M_{\kappa}$--polyhedral complex with only finitely many isometry types of cells.  If $\kappa \leq 0$, then $X$ is $\CAT (\kappa)$ if and only if $X$ satisfies the link condition and is simply connected.  If $\kappa =1$, then $X$ is $\CAT(1)$ if and only if $X$ satisfies the link condition and contains no isometrically embedded circles of length less than $2\pi$.
\end{theorem}

If $\kappa \leq 0$, to prove that an $M_{\kappa}$--complex is $\CAT (0)$, it thus suffices to show that it is simply connected and has $\CAT(1)$ links.  Applying \fullref{linkcond} to the links then reduces the problem to one of ruling out short loops in successive links.  In dimension 2, links of vertices are graphs, where checking for embedded loops is straightforward.

\begin{prop}{\rm\cite[Lemma II.5.6]{BridsonHaefliger}}\qua
A $2$--dimensional $M_{\kappa}$--complex satisfies the link condition if and only if every injective loop in the link of every vertex of $X$ has length at least $2 \pi$. \end{prop}

\subsection{Isometries and group actions}
\begin{defn} \label{def:elliptic} By an {\em action} of a group $G$ on a space $X$, we mean a homomorphism $\phi \co G \rightarrow \Isom(X)$.  (We therefore consider only isometric actions, but we will not, in general, assume that actions are proper or faithful.)  A group action on $X$ is {\em trivial} if a point of $X$ is fixed by every element of the group.  An isometry $\gamma \in \Isom (X)$ is {\em semisimple} if $d(\cdot, \gamma (\cdot))$ attains a minimum on $X$ and an isometry is {\em elliptic} if it fixes a point.  For $H\subset G$, we say $H$ is {\em $\phi$--elliptic} if the action of $H$ on $X$ given by $\phi \mid _H$ is trivial.
\end{defn}

The following lemma consists of standard facts about fixed sets.
\begin{lemma}\label{fixedsetfacts}
Suppose a group $G$ acts by isometries on a geodesic metric space $(X,d)$ via $\phi \co G \rightarrow \Isom(X)$.
\begin{enumerate}
\item[(i)] If $d$ is convex (eg if $X$ is $\CAT(0)$), then for any $\phi$--elliptic subset $H \subset G$, the set of global fixed points of $H$ under $\phi$, denoted $\Fix _{\phi} (H)$, is contractible.
\item[(ii)]  If $H_1 , H_2 \subset G$, then $\Fix _{\phi} ( \langle H_1 \cup H_2 \rangle ) = \Fix_{\phi}(H_1) \cap \Fix_{\phi}(H_2)$.
\end{enumerate}
\end{lemma}

The next result \cite[Corollary II.2.8]{BridsonHaefliger} is crucial for our arguments.

\begin{prop}[Bruhat--Tits Fixed Point Theorem] \label{BT}
Let $X$ be a complete, connected $\CAT (0)$ space.  Any action of a finite group on $X$ is trivial.  More generally, any action on $X$ with a bounded orbit is trivial.
\end{prop}

The proof of this result depends on the fact that the center of a bounded set is well-defined in complete $\CAT (0)$ spaces.  Since every orbit is preserved, the center of a bounded orbit, which is metrically defined, is preserved as well.

A standard consequence is the following:
\begin{cor}\label{finiteindexenough}
Suppose a group $G$ acts by isometries on a complete, connected $\CAT (0)$ space $X$ via $\phi \co G \rightarrow \Isom (X)$. If $H$ is a finite index subgroup of $G$ and $H$ is $\phi$--elliptic, then $G$ is $\phi$--elliptic.
\end{cor}
\begin{proof}
Since $H$ is finite index in $G$, there is a normal subgroup $N$ of finite index in $G$ such that $N<H<G$.  Now $N$ is $\phi$--elliptic since $H$ is.  Denote the (nonempty, contractible) fixed set of $N$ by $X^N$.  Then $G/N$ acts on $X^N$.  This action has a fixed point $x_G$ by \fullref{BT}.  The point $x_G$ is fixed by the group $G$.
\end{proof}

\subsection{Actions of Coxeter groups} \label{Coxactions}
Coxeter groups act on $\CAT (0)$ spaces.  In this section, we consider some of these actions.  Fix a Coxeter system $(W,S)$ and suppose $W$ is infinite.
\begin{itemize}
\item The {\em Davis--Moussong complex} $\Sigma_{DM}$ associated to $(W,S)$ (see, for example, \cite{DavisMoussong,Davisbook}) has a simplicial structure similar to that of the classical Coxeter complex.  The dimension of $\Sigma_{DM}$ is given by the maximal rank of finite special subgroups of $W$.  In his thesis \cite{Moussong}, Moussong proved that by assigning appropriate Euclidean metrics to the cells of a particular cellular structure on $\Sigma_{DM}$, the resulting metric on $\Sigma_{DM}$ is $\CAT (0)$.  The group $W$ acts on $\Sigma_{DM}$ properly by isometries.
\item Given another Coxeter system $(W',S)$, there is a surjection of $W$ to $W'$ if $m'_{ij}$ divides $m_{ij}$ for all $i,j$.  In particular, $W$ acts nontrivially on the Davis--Moussong complex of each such infinite quotient group.
\item Niblo and Reeves construct in \cite{NibloReeves} a locally finite, finite-dimensional $\CAT (0)$ cube complex on which $W$ acts properly discontinuously.  For right-angled Coxeter groups, their complex is isometric to the Davis--Moussong complex.
\item Since $W$ is generated by torsion elements, it has no nontrivial homomorphisms to $\ZZ$.  However, Cooper, Long and Reid in \cite{CLR} and Gonciulea in \cite{Gonciulea} showed that $W$ has a finite index subgroup that surjects to $\ZZ$.  Such a finite index subgroup then acts nontrivially by translation on the real line.  However, $W$ itself acts nontrivially on a tree if and only if some $m_{ij}$ is infinite (see \fullref{FACox} below).
\end{itemize}

\section{Group decompositions}
\subsection{Amalgamated products and normal form length} \label{amalgcox}
Recall from Bass--Serre Theory (see \cite{scottwall} or \cite{SerreTrees}) that actions on trees correspond to decompositions of groups as graphs of groups.  We denote by $A*_CB$ the {\em amalgamated product} of $A$ and $B$ along $C$.  Thus $A*_CB$ is the pushout of the diagram of groups $A\hookleftarrow C \hookrightarrow B$.  That is, it is the fundamental group of the graph of groups with vertices $A$ and $B$ and a single edge $C$.

\begin{defn}
Let $G=A*_CB$, and let $g \in G$.  The {\em normal form length} of $g$ is defined to be the following:
\[ \min \{ k \st g =\underbrace{a_1b_2a_3\ldots}_k\textrm{ or }\underbrace{b_1a_2b_3\ldots}_k
\textrm{ with }a_i\in A \textrm{ and } b_j \in B\  \forall \ 1\leq i,j \leq k\}.\]  Note that this is equal to the length of the normal form representative of $g$ in $A*_CB$ (see \cite{scottwall} or \cite{SerreTrees} for details on normal forms in amalgamated products).  We will be particularly interested in the normal form length of words in amalgamated products of Coxeter groups.
\end{defn}

Let $(W,S)$ be a Coxeter system with Coxeter matrix $M$ and special subgroups $A$ and $B$.  Let $M_{AB}$ denote the submatrix of $M$ corresponding to the generators $S_A\cup S_B$.  Define $M'$ to be the matrix $M_{AB}$ with only the following change:  for $s_a \in S_A - S_B$ and $s_b \in S_B -S_A$, the corresponding entry in $M'$ is $ \infty$. Then $(A*_{A\cap B}B, S_A\cup S_B)$ is a Coxeter system with Coxeter matrix $M'$.

Let $w=s_1s_2\cdots s_k$ be a word in the alphabet $S_A \cup S_B$.  We
denote by {\em $\alt(w)$} the minimum number $k$ such that there are
indices $1= i_1 <i_2< \cdots <i_k\leq r$ with
$s_{i_j},s_{i_j+1},\ldots ,s_{i_{j+1}-1} \in S_A$ for $j$ odd and
$s_{i_j},s_{i_j+1},\ldots ,s_{i_{j+1}-1} \in S_B$ for $j$ even (or
vice versa).  That is, $\alt (w)$ is the number of ``alternations''
between elements of $A$ and elements in $B$ in the word $w$.  The
normal form length of an element $g\in G$ is thus $\min \{ \alt(w) \st
w\mbox{ is a word representing }g\}$.

Note that the only Type (II) elementary $M'$--operations in $G$ are those in $A$ and $B$ themselves.  Together with \fullref{titswp}, this implies the following:

\begin{lemma}\label{coxred} \label{nflred} \label{lemma:alt}
Let $W$ be a Coxeter group with special subgroups $A$ and $B$, and let $G=A*_{A\cap B}B$.  Let $M'$ denote the Coxeter matrix of $(G,S_A\cup S_B)$ as defined above.
\begin{enumerate}
\item[(i)] Let $w$ be a word in the alphabet $S_A\cup S_B$ and let $w'$ be any subword of a word obtained from $w$ via a sequence of elementary $M'$--operations.  Then $\alt (w')\leq \alt(w)$.
\item[(ii)] For $g\in G$, every $M'$--reduced word representing $g$ realizes the normal form length of $g$.  That is, if the normal form length of $g$ in $G$ is $k$ then $\alt(w)=k$ for every $M'$--reduced word $w$ representing $g$.
\item[(iii)] For $a\in S_A-S_B$ and $b\in S_B - S_A$, the normal form length in $G$ of an alternating product of $a$ and $b$ is the length of the product as a word in the alphabet  $\{a,b\}$.
\end{enumerate}
\end{lemma}
\vspace{3pt}

\subsection{Triangles of groups}\label{GerStsection}
\vspace{3pt}
The study of triangles of groups is a $2$--dimensional analogue of Bass--Serre Theory due originally to S Gersten and J Stallings (see \cite{StallingsTriangles}).  In the next section, we will consider a more general construction, but the case of triangles will be of particular interest.
\vspace{3pt}

Suppose a group $G$ acts on a $2$--dimensional simplicial complex with quotient a triangle.  By choosing a fundamental domain for the action, we can assign stabilizer groups to the vertices, edges, and face of this triangle.  This process yields a triangle of groups, a $2$--dimensional analogue of the graph of groups corresponding to $A*_CB$.  Formally, a {\em triangle of groups} is a commutative diagram of groups and monomorphisms, as in \fullref{triangleofgroups} below.  The groups $A$, $B$ and $C$ are {\em vertex groups}, $D$, $E$ and $F$ are {\em edge groups}, and $K$ is the {\em face group}.  The vertex, edge and face groups are all called {\em local groups}.

\begin{figure}[ht!]
\begin{center}
\includegraphics{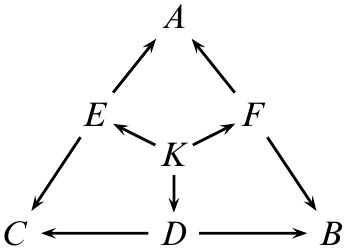}
\caption{Triangle of Groups}\label{triangleofgroups}
\end{center}
\end{figure}
\vspace{3pt}

The {\em fundamental group} $G$ of the triangle of groups is the colimit of the diagram, ie, the unique group (up to isomorphism) satisfying the following universal mapping property:
\begin{quote}
Given any group $H$ and homomorphisms from the vertex groups to $H$, there is a unique homomorphism from $G$ to $H$ such that the resulting diagram of groups commutes.
\end{quote}
\vspace{3pt}

\begin{defn}
Let $A$ and $K$ be groups, $E$ and $F$ be subgroups of $A$, and $K\rightarrow E$ and $K\rightarrow F$ be homomorphisms. Let $\phi$ be the natural surjection $\phi\co E*_KF\rightarrow \langle E,F \rangle \subset A$.  The {\em Gersten--Stallings angle between $E$ and $F$ over $K$}, denoted $\angle_A (E,F;K)$, is $\frac{2\pi}{n}$ where \[ n=\min \{\mbox{normal form lengths of nontrivial elements in }\ker(\phi)\}.\]  Equivalently, $n$ is the length of the shortest loop in the graph whose vertices are the cosets of $E$ and $F$ in $A$ and whose edges are the cosets of $K$ in $A$, with incidence given by inclusion.  Note that $n$ is {\em even} since this coset graph is bipartite.
In the triangle of groups of \fullref{triangleofgroups}, we say that the {\em Gersten--Stallings angle at the vertex $A$} is $\angle_A (E,F;K)$.
\end{defn}

Given a triangle of groups with fundamental group $G$ such that all the vertex groups inject into $G$, there is a natural $2$--complex $X$, the {\em universal cover} of the triangle of groups, on which $G$ acts with quotient a triangle.  The $2$--simplices of $X$ are given by cosets of $K$ in $G$, with incident edges given by the corresponding cosets of $D$, $E$ and $F$, and vertices the cosets of $A$, $B$ and $C$.  We may give $X$ a metric by assigning a metric on the fundamental domain that agrees with the Gersten--Stallings angles.  In particular, each triangle can be given the metric of a triangle in $\EE ^2$, $\HH ^2$ or $\SS ^2$, according to the sum of the Gersten--Stallings angles at each vertex of the triangle of groups.  The resulting metric on $X$ is then piecewise-Euclidean, piecewise-hyperbolic or piecewise-spherical.

The assigned metric ensures that $X$ satisfies the link condition (see \fullref{linkconddefn}).  This implies the following theorem of Gersten and Stallings (see \cite{StallingsTriangles}).

\begin{theorem} \label{CAT0triangles}
If the sum of the Gersten--Stallings angles in a triangle of groups is less than or equal to $\pi$, then its universal cover $X$ is $\CAT (0)$.  If the sum is strictly less than $\pi$, then $X$ is $\CAT(-1)$.
\end{theorem}

\subsection{Simple complexes of groups}\label{scog}
The theory of {\em complexes of groups} was studied by Haefliger in \cite{Haefliger} and independently (in dimension 2) by Corson in \cite{corson}.  We will follow the development in \cite{BridsonHaefliger} but will only need the theory of simple complexes of groups.

\begin{defn}
A {\em simple complex of groups} over a poset $\mathcal Q$, denoted $G(\mathcal Q)$, consists of:
\begin{enumerate}
\item For each $\sigma \in {\mathcal Q}$, a group $G_{\sigma}$ (the {\em local group at $\sigma$}).
\item For each $\tau < \sigma$, a monomorphism $\iota_{\tau \sigma} \co G_{\sigma} \rightarrow G_{\tau}$ so that   $\iota _{\tau \sigma} = \iota _{\tau \rho}  \iota _{\rho \sigma}$ whenever $\tau < \rho < \sigma$.
\end{enumerate}

Denote by $\widehat{G({\mathcal Q})}$ the colimit of the diagram of groups $G({\mathcal Q})$.  We will take $\mathcal Q$ to be the poset associated to the cells of a cell complex.  Note that if the geometric realization of $\mathcal Q$ is simply connected, as is true in the cases we will consider, then $\widehat{G({\mathcal Q})}$ is also what is known as the {\em fundamental group} of the corresponding complex of groups.
\end{defn}

\begin{remark} \label{prescofg}
Given presentations for the local groups of $G({\mathcal Q})$, we have a presentation for the colimit $\widehat{G ( {\mathcal Q})}$.  The generators are given by the generators of the local groups, and along with the relations coming from the local groups, elements of different local groups are related via the monomorphisms $\iota _{\tau \sigma}$.  Specifically, if the local groups have presentations $G_{\sigma} = \langle S_{\sigma} \mid R_{\sigma}\rangle$, then we have \[\widehat{G({\mathcal Q})}= \left< \bigcup_{\sigma \in \mathcal Q} S_{\sigma} \biggm| \bigcup_{\sigma \in \mathcal Q} \left\{ R_{\sigma} , \{g=\iota_{\tau \sigma} (g) \st g \in G_\sigma ,\ \tau < \sigma \} \right\} \right>.\]
\end{remark}

Suppose $\mathcal Q$ is the poset corresponding to the simplices of a $2$--simplex, ordered by inclusion.  A simple complex of groups over $\mathcal Q$ is a triangle of groups.  More generally, if $\mathcal Q$ is the poset corresponding to the simplices of an $n$--simplex $\Delta ^n$, ordered by inclusion, we refer to a simple complex of groups over $\mathcal Q$ as an {\em $n$--simplex of groups}.

As was the case for triangles, if the local groups inject into $\widehat{G({\mathcal Q})}$, there is a simply connected space, the {\em universal cover of $G({\mathcal Q})$}, upon which $\widehat{G({\mathcal Q})}$ acts with quotient the underlying cell complex.  In particular, the underlying cell complex is a {\em strict} fundamental domain for the action.  Assigning a metric to the underlying cell complex of $\mathcal Q$ yields a metric on the universal cover and the resulting action is by isometries.  We will call a simplex of groups {\em nonpositively curved} if its universal cover is $\CAT (0)$.

The following is an analogue of the graph of groups version given by Mihalik and Tschantz in \cite{MihalikTschantz}.

\begin{prop} \label{complexofcoxeterprop}
Let $\mathcal Q$ be the poset of cells of a complex $\overline{\mathcal Q}$ ordered by inclusion.  Suppose $(W,S)$ is a Coxeter system and $G(\mathcal Q)$ is a simple complex of special subgroups of $W$ with monomorphisms $\iota _{\tau \sigma}$ given by natural inclusions.  Then, $\widehat{G(\mathcal Q)} = W$ (with all the resulting homomorphisms from the local groups to $W$ the natural inclusions) if both of the following hold:
\begin{enumerate}
\item[(i) ] For $s \in S$, the set $\{ \sigma \st \sigma \in \mathcal Q, \ s \in G_{\sigma}\}$ corresponds to a nonempty connected subcomplex of $\overline{\mathcal Q}.$
\item[(ii)] If $s,t \in S$ and $m_{st}<\infty$, then $\{s,t\} \subset G_{\sigma}$ for some $\sigma \in \mathcal Q$.
\end{enumerate}
\end{prop}

\begin{proof}
For each local group $G_{\sigma}$, let $\langle S_{\sigma} \mid R_{\sigma}\rangle$ be the induced Coxeter presentation of $G_{\sigma}$ as a special subgroup of $W$.

By \fullref{prescofg}, we have \[\widehat{G({\mathcal Q})}= \left< \bigcup_{\sigma \in \mathcal Q} S_{\sigma} \biggm| \bigcup_{\sigma \in \mathcal Q} \left\{ R_{\sigma} , \{g=\iota_{\tau \sigma} (g) \st g \in G_\sigma ,\ \tau < \sigma \} \right\} \right>.\]

Suppose that (i) and (ii) hold.  Since the monomorphisms of $G(\mathcal Q)$ are the natural inclusions, by (i) we can write
\[\widehat{G({\mathcal Q})} = \left< S \biggm| \bigcup_{\sigma \textrm{ vertex of }\overline{\mathcal Q}} R_{\sigma} \right>.\]
So, it suffices to show that the relations appearing in the vertex groups are precisely the relations of $W$.  By our choice of presentations for the local groups (as special subgroups of $W$), every relation $r_{\sigma} \in R_{\sigma}$ is a relation from our Coxeter presentation for $(W,S)$.  Moreover, by (i), every relation $s^2=1$ appears in some local group (hence some vertex group).  Finally, by (ii), every relation of the form $(st)^{m_{st}} =1$ appears in some local group (so in a vertex group).
\end{proof}

\vspace{4pt}
\section[Property FA _n]{Property $\FA _n$}
\vspace{4pt}
\subsection{Definition and examples}
\vspace{4pt}

A group $G$ has Serre's {\em  property $\FA$} (see \cite{SerreTrees}) if for every tree $T$, every action (without inversions) of $G$ on $T$ is trivial (has a global fixed point).  In particular, if $G$ has FA, then $G$ does not split nontrivially as a graph of groups.  Since trees are $\CAT (0)$, the Bruhat--Tits Theorem (\fullref{BT}) applies, so finite groups have FA.  Other groups known to have $\FA$ include the following:
\begin{itemize}
\item Every finite index subgroup of SL$_n(\ZZ)$ for $n \geq 3$ (Margulis--Tits, see \cite{SerreTrees})
\item Finitely generated torsion groups (Serre \cite{SerreTrees})
\item Coxeter groups such that $m_{ij} < \infty$ for all $i$ and $j$ (Serre \cite{SerreTrees})
\item Finitely generated groups with Kazhdan's property (T) (Watatani \cite{watatani})
\item Out$ (F_n)$ and Aut$(F_n)$ for $n \geq 3$ (Bogopolski \cite{Bogopolski}, Culler--Vogtmann \cite{CullerVogtmann})
\item Mapping class groups of higher genus surfaces (Culler--Vogtmann \cite{CullerVogtmann})
\end{itemize}

\begin{remark} \label{FACox}
Mihalik and Tschantz note in \cite{MihalikTschantz} that a Coxeter group has FA if and only if every $m_{ij}$ is finite.  In particular, if $(W,S)$ is a Coxeter system and $m_{ij} = \infty$, then \hbox{$W \cong \langle S-\{s_i\} \rangle *_{\langle S - \{s_i, s_j\}\rangle} \langle S- \{s_j\} \rangle$.}
\end{remark}

\begin{defn}  A {\em $\CAT(0)$ $n$--complex} is an $M_{\kappa}$--polyhedral $n$--complex (see \fullref{CAT0basics}) that is complete, connected, $\CAT (0)$, and has only finitely many isometry types of cells.  A group $G$ has {\em property $\FA _n$} if for every $\CAT (0)$ \hbox{$n$--complex} $X$, every action of $G$ on $X$ by cellular isometries has a global fixed point.  The group $G$ has {\em strong $\FA _n$} if for every complete, connected $\CAT (0)$ space $X$ of topological dimension $n$, every action of $G$ on $X$ by semisimple isometries has a global fixed point.
\end{defn}

Bridson showed in \cite{Bridson99} that if $X$ is a connected $M_{\kappa}$--polyhedral complex with only finitely many isometry types of cells, then every cellular isometry of $X$ is semisimple.  Therefore, any group with strong $\FA_n$ also has $\FA_n$.

The following, as noted by Farb in \cite{FarbHelly}, is an immediate consequence of \fullref{BT}.

\begin{cor}\label{finiteFAn}
Let $G$ be a finite group.  Then $G$ has strong $\FA_n$ for all $n$.
\end{cor}
\vbox{
Farb also showed in \cite{FarbHelly} that the following groups have property $\FA _n$:
\begin{itemize}
\item Finite index subgroups of SL$_m ( \ZZ)$ and SL$_m(\ZZ [1/p])$ for $m\geq n+2$
\item More generally, arithmetic or $S$-arithmetic subgroups of $K$--simple algebraic $K$--groups of $K$--rank at least $n+1$ for $K$ a global field
\item Discrete groups generated by reflections in the sides of Euclidean or hyperbolic $(n+1)$--simplices
\end{itemize}}
Recently, Bridson \cite{BridsonHelly} has also studied property $\FA_n$ for automorphism groups of free groups.

\subsection{Homological techniques and implications} \label{homimp}
Note that throughout we consider homology with $\ZZ$ coefficients.

We will prove that certain actions are trivial by studying the combinatorics of fixed sets.  Consider a collection of sets $\{S_{\alpha} \}_{\alpha \in I}$.  Recall that the {\em nerve} of this collection, denoted $\cN (\{S_{\alpha}\}_{{\alpha}\in I})$,  is the simplicial complex whose vertices are indexed by the set $I$ and such that the set of vertices corresponding to $J\subset I$ span a simplex if and only if $ \bigcap_{\alpha \in J} S_{\alpha } \neq \emptyset$.

\begin{notation}  We will be particularly interested in the nerves of collections of fixed sets.  Let $G$ be a group, and suppose $\phi\co G \rightarrow \Isom(X)$ is the homomorphism describing an action of $G$ on $X$.  Let $\cS$ be a finite collection of $\phi$--elliptic subsets of $G$ (see \fullref{def:elliptic}).  We denote by $\cN(\cS, \phi)$ the nerve of the collection of fixed sets of the elements of $\cS$.  That is $\cN(\cS, \phi) = \cN\left(\{\Fix_{\phi}(S)\}_{S \in \cS}\right).$
\end{notation}

\begin{remark} \label{simplextrivial}
Suppose $\cS$ is finite.  Then $\cN (\cS , \phi)$ is a simplex if and only if there is a global fixed point for the action the subgroup $\langle \Sigma \rangle$ of $G$ generated by the union of the elements of $\cS$.  By \fullref{finiteindexenough}, if $\langle \cS \rangle$ has finite index in $G$, then this holds if and only if the action of $G$ has a global fixed point.
\end{remark}

Our main technique is based on the following two results.

\begin{theorem}{\rm(Leray \cite[Theorem VII.4.4]{BrownCohom})}\label{CWleray}\qua
Suppose $X$ is a CW complex and is the union of subcomplexes $X_{\alpha}$ such that the intersection of any finite subcollection of the $X_{\alpha}$ is either empty or acyclic.  Then $H_*(X) = H_* \left(\cN(\{X_{\alpha}\})\right)$.
\end{theorem}

\begin{theorem}{\rm(McCord  \cite[Theorem 2]{mccord})}\qua
\label{leray}%
Let $X$ be a space and $\mathcal U$ a locally finite open cover of $X$ such that the intersection of any finite subcollection of $\mathcal U$ is either empty or homotopically trivial.  Then, there is a weak homotopy equivalence $\cN(\mathcal U) \rightarrow X$ so $H_*(\cN(\mathcal U))=H_*(X)$.
\end{theorem}

\begin{defn}
Let $K$ be a simplicial complex.  We say that $K$ is {\em n--allowable} if  $H_m(K) = 0$ for all $m \geq n$.
\end{defn}

The motivation for this definition can be found in the following, which is implicit in \cite{FarbHelly}.  We include a proof here for completeness.

\begin{prop} \label{leraygroupcriterion}
Let $G$ be a group.  Suppose $G$ acts on a complete $\CAT (0)$ space $X$ of dimension $n$ via $\phi\colon G \rightarrow \Isom(X)$.  Let $\cS$ be a finite set of $\phi$--elliptic subsets of $G$.  Then $\cN (\cT , \phi)$ is $n$--allowable for all $\cT \subset \cS$.
\end{prop}

\begin{proof}
Let $Y$ be the union of the fixed sets $\Fix _{\phi} (S)$ for $S\in \cT$.  By \fullref{fixedsetfacts}, nonempty intersections of the sets $\Fix_{\phi}(S)$ are also fixed sets so are contractible.  Taking regular neighborhoods of the sets $\Fix _{\phi}(S)$ that preserve the intersection data, we apply \fullref{leray} to the resulting open cover to find that \hbox{$H_* (\cN (\cT , \phi )) = H_* (Y)$.}

Let $m\geq 1$.  Applying the long exact homology sequence for pairs, we have the exact sequence $H_{m+1}(X) \rightarrow H_{m+1} (X,Y) \rightarrow H_m (Y) \rightarrow H_m(X).$  Since $X$ is contractible (see \fullref{CAT0props}), we find that $H_{m+1} (X,Y) \cong H_m (Y)$.   Now  $H_{m+1} (X,Y) = 0$ for $m\geq n$ as $X$ is $n$--dimensional.  Thus for $m\geq n$ we have $H_m(Y) =0$ so $H_m ( \cN ( \cT , \phi))=0$ as well.
\end{proof}

Note that in the above proof we only used the $\CAT (0)$ assumption for contractibility of fixed sets, so the result holds more generally.

\begin{remark} \label{simplexleray}
Let $K$ be an $n$--allowable simplicial complex with $0$--skeleton consisting of the vertices $v_0, \ldots, v_{k}$ for some $k \geq n$.  Then $H_k(K)=0$ by the definition of $n$--allowable.  So, if the $(k-1)$--skeleton of $K$ is the boundary of a $k$--simplex, then the $k$--skeleton (and hence the entire complex $K$) is a $k$--simplex.
\end{remark}

Using a topological form of Helly's Theorem, Farb applies certain cases of the following result in \cite{FarbHelly}.   We include here a proof of the general result using the techniques from above.
\begin{cor} \label{stronghelly}
Let $G$ be a group and $\cS$ a finite collection of subsets of $G$ whose union generates a finite index subgroup of $G$.  Suppose $G$ acts on a complete $n$--dimensional $\CAT (0)$ space $X$ via $\phi \co G \rightarrow \Isom(X)$.  If $n<| \cS |$ and if every $n+1$ elements of $\cS$ generate a $\phi$--elliptic subgroup of $G$, then $G$ is $\phi$--elliptic.
\end{cor}

\begin{proof}
By assumption, every $n+1$ elements of $\cS$ generate a $\phi$--elliptic subgroup.  So, by definition of $\cN = \cN (\cS , \phi)$, every $n+1$ vertices in $\cN$ span an $n$--simplex.  If $|\cS|=n+1$, then $G$ is $\phi$--elliptic by \fullref{simplextrivial}.  Otherwise, consider a subset $\cT \subset \cS$ of cardinality $n+2$.  Then, the $n$--skeleton of $\cN(\cT , \phi)$ is the boundary of an $(n+1)$--simplex.  By \fullref{leraygroupcriterion}, we know that $\cN(\cT , \phi)$ is $n$--allowable, so by \fullref{simplexleray}, its $(n+1)$--skeleton is actually an $(n+1)$--simplex.  Inductively, we see that the $(|\cS|-1)$--skeleton of $\cN$ is an $(|\cS|-1)$--simplex, so $\cN$ is a simplex.  Thus $G$ is $\phi$--elliptic by \fullref{simplextrivial}.
\end{proof}

As an immediate consequence we have the following:

\begin{cor}\label{grouphelly}
\label{thmgenhelly}
Let $G$ be a group and let $\cS$ be a finite collection of subsets of $G$ whose union generates a finite index subgroup of $G$.  If for every subcollection $\{S_0, S_1, \ldots , S_n\} \subset \cS$ there is an $\FA _n$ (resp.\ strong $\FA _n$) subgroup $H\subset G$ such that $S_i \subset H$ for $0\leq i\leq n$, then $G$ has $\FA _n$ (resp.\ strong $\FA _n$).  In particular, if every $n+1$ elements of $\cS$ generate a group with property $\FA _n$ (resp.\ strong $\FA _n$), then $G$ has property $\FA _n$ (resp.\ strong $\FA _n$).
\end{cor}

\begin{defn}
(Notation as in \fullref{scog}.)  We say that an $n$--simplex of groups is {\em minimal} if for all $k<n$, every local group $G_{\sigma}$ corresponding to a $k$--simplex $\sigma$ is generated by the (images of the) local groups $G_{\tau}$ such that $\sigma \subsetneq {\tau}$.  We will call such a $G_{\sigma}$ a {\em local group of codimension $n-k$}.
\end{defn}
A version of the following was proved for $\FA _1$ by R Alperin in \cite{alperin}.

\begin{cor}\label{minimalsimplex}
Let $n\geq 1$.  Suppose a group $G$ has a decomposition as a minimal $(n+1)$--simplex of groups $\Lambda$ such that every local group of $\Lambda$ has $\FA_n$ (resp.\ strong $\FA _n$).  Then, $G$ has $\FA_n$ (resp.\ strong $\FA _n$).
\end{cor}

\begin{remark} Note that the simplex of groups decomposition does not have to be nonpositively curved.  This result implies, for example, that the fundamental group of any realizable minimal $n$--simplex of finite groups has strong $\FA _{n-1}$.  So, if a group acts nontrivially by semisimple isometries on any complete, connected $\CAT (0)$ space of dimension $n$, then it does not decompose as a minimal $m$--simplex of finite groups for any $m>n$.
\end{remark}

\begin{proof}[Proof of \fullref{minimalsimplex}]
Let $ \cS = \left\{\mbox{local groups of codimension 1}\right\}.$  By assumption, $H$ has $\FA_n$ for all $H \in \cS$ and by minimality, the union of the elements of $\cS$ generates $G$.  Moreover, for all $k\leq n+1$, every collection of $k$ elements of $\cS$ generates a local group of codimension $k$, which in turn has $\FA_n$.  Thus $G$ has $\FA_n$ by \fullref{grouphelly}.  The case of strong $\FA _n$ is analogous.
\end{proof}

\section[Strong FA _n for Coxeter groups]{Strong $\FA _n$ for Coxeter groups}
\subsection{$n$--spherical Coxeter groups}
\begin{defn}
Let $(W,S)$ be a Coxeter system.  The Coxeter group $W$ is {\em $n$--spherical} if every special subgroup of $W$ of rank less than or equal to $n$ is a finite group.  (Note that for $n\leq |S|$, the group $W$ is $n$--spherical if and only if all of its special subgroups of rank $n$ are finite.)  A special subgroup $W_T$ is {\em $n$--spherical} if it is $n$--spherical with respect to the generating set $T$.
\end{defn}

Finite groups have $\FA_n$ by \fullref{finiteFAn}, so \fullref{thmhelly} follows from \fullref{grouphelly} by setting $\cS = \{ \{s\} \st s\in S\}$.  In fact, \fullref{stronghelly} implies the following:
\begin{cor} \label{Coxnfinite}
Let $W$ be a Coxeter group.  If $W$ is $(n+1)$--spherical, then every action of $W$ on every complete connected $\CAT (0)$ space of dimension $n$ has a global fixed point.  In particular, $W$ has strong $\FA _n$.
\end{cor}

The converse is true in dimension 1 (see \cite{MihalikTschantz}), and we will prove it also holds in dimension 2.  For higher dimensions, we have the following reformulation of \fullref{mainconj}.

\begin{conj}[Coxeter $\FA _n$ Conjecture] \label{pdimconj}
Let $(W,S)$ be a Coxeter system.  The following are equivalent:
\begin{enumerate}
\item[(i)] $W$ is $(n+1)$--spherical.
\item[(ii)] Every action of $W$ on every complete connected $\CAT (0)$ space of dimension $n$ has a global fixed point.
\item[(iii)]  $W$ has strong $\FA_n$.
\item[(iv)] $W$ has property $\FA_n$.
\item[(v)]  $W$ does not split nontrivially as a nonpositively curved $m$--simplex of special subgroups, for all $0< m \leq n$.
\end{enumerate}
\end{conj}

We have already seen that (i) $\implies$ (ii).  That (ii) $\implies$
(iii) $\implies$ (iv) $\implies$ (v) is clear by definition.  It
remains only to show that (v) $\implies$ (i).  We describe an approach
to this in the next section.

Note that property $\FA_n$ is, by definition, a property of a group rather than of a presentation of the group.  However, by the classification of Coxeter groups, the property of being $n$--spherical is easily verified by looking at a Coxeter presentation (see \fullref{lookupfinite}).  Thus, if the Coxeter $\FA _n$ Conjecture holds in general, then property $\FA _n$ can be detected by considering any Coxeter presentation.

\subsection[The CAT(0) Conjecture]{The $\CAT (0)$ Conjecture} \label{construct}
We now describe an approach to completing the proof of the Coxeter $\FA _n$ Conjecture (\fullref{pdimconj}).  In particular, the goal is to prove that if $W$ is not $(n+1)$--spherical, then $W$ splits nontrivially as a nonpositively curved simplex of special subgroups of dimension at most $n$.  Note that if $W$ is finite, then $W$ is $(n+1)$--spherical for all $n$.  We thus only need consider infinite Coxeter groups in what follows.

Fix a Coxeter system $(W,S)$ with $W$ infinite.  Let \[v=v(W)= \max\{m \st W \textrm{ is $m$--spherical}\}.\]  Our goal is to construct a decomposition of $W$ as a $v$--simplex of groups.  To do this, we first consider a ``natural'' decomposition of a subgroup of $W$.

Let $S' \subset S$ be the generating set of an infinite special subgroup of $W$ of rank $v+1$.  Note that such a subset exists by the definition of $v$.  Moreover, every proper subset of $S'$ generates a finite group.  By \fullref{classification}, the Coxeter group $W'$ is thus an irreducible Euclidean or hyperbolic reflection group with fundamental domain a compact simplex.  Hence, $W'$ has a natural decomposition as a Euclidean or hyperbolic $v$--simplex of groups $\Lambda '$ in which the vertex groups are the special subgroups of $W'$ of rank $v$.  We assign to the simplex $\Lambda '$ the metric given by the Euclidean or hyperbolic metric on a fundamental domain for the natural action of $W'$.  Then, the universal cover $X'$ of $\Lambda '$ is isometric to the original hyperbolic or Euclidean space, and the action of $W'$ on $X'$ is the original action.

Following the motivation from dimension one (see \fullref{FACox}), we construct a new simplex of groups $\Lambda$ by adding the missing generators.  In particular, $\Lambda$ is given by the same metric simplex as was $\Lambda '$, but for $A\varsubsetneq S'$, the local group $\langle A \rangle$ is replaced by the group $\langle{ A \cup (S-S')}\rangle$.  Note that the local group associated to the maximal simplex in $\Lambda $ is thus $\langle{S-S'}\rangle$.

\begin{eg}
Let $S=\{a,b,c,x,y,z\}$ and let $W$ be the Coxeter group given by the Coxeter diagram in \fullref{2dimdiagram}.

\begin{figure}[ht!]
\begin{center}
\includegraphics{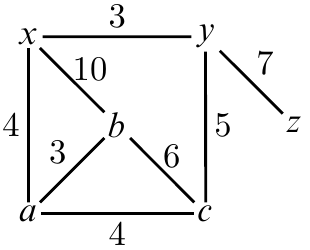}
\caption{Coxeter diagram with $v(W)=2$}\label{2dimdiagram}
\end{center}
\end{figure}

Then $W$ is $2$--spherical since every $m_{ij}$ is finite.  Let $S' = \{ a, b, c \}$.  Then $W' = \langle S' \rangle$ is infinite, so $v(W) = 2$.  Moreover, $W'$ is a hyperbolic triangle group with angles $\frac{\pi}{3}, \frac{\pi}{4}$ and $ \frac{\pi}{6}$, and the action of $W'$ on $\HH ^2$ gives the triangle of groups $\Lambda '$ shown in \fullref{hyptridecomp}.  The resulting triangle of groups $\Lambda$ is shown in \hyperlink{double1}{Figure 5}.  The metric on $\Lambda$ is that of a hyperbolic $(3,4,6)$--triangle.\end{eg}

\begin{figure}[ht!]
\begin{minipage}{0.47\linewidth}
\hypertarget{double1}{\cl{\includegraphics{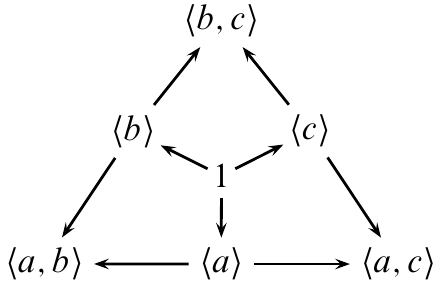}}}
\caption{Natural\newline $2$--splitting $\Lambda '$ of $W'$}\label{hyptridecomp}
\end{minipage}
\addtocounter{figure}{1}
\begin{minipage}{0.5\linewidth}
\cl{\includegraphics{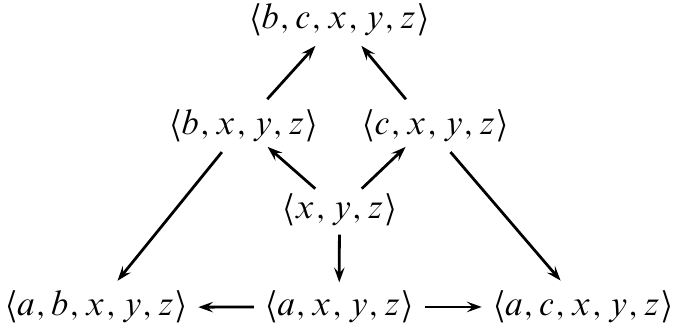}}
\caption{$2$--splitting $\Lambda$\newline determined by $S'$}\label{2dimegdecomp}
\end{minipage}
\end{figure}

As an immediate consequence of \fullref{complexofcoxeterprop} we have the following:
\begin{cor}
Let $\Lambda$ be the simplex of groups associated to a rank $v+1$ infinite special subgroup of $W$ as constructed above.  Then  $\Lambda$ is a splitting of $W$.  That is $W=\widehat{\Lambda}$ (notation as in \fullref{scog}).
\end{cor}

We will refer to the decomposition $\Lambda$ of $W$ as the {\em v--splitting of W determined by $S'$}.  The group $W$ acts on the simply connected universal cover of $\Lambda$ with quotient a $v$--simplex.

\begin{CATOconj} \label{cat0conj}
In the above construction, the universal cover of the $v$--splitting of $W$ determined by $S'$ is $\CAT (0)$.
\end{CATOconj}

The preceding discussion actually proves the following:

\begin{theorem}\label{CAT0impliesFAn}
$\CAT (0)$ Conjecture $\implies$ Coxeter $\FA _n$ Conjecture (\fullref{pdimconj}).
\end{theorem}

\begin{remark}\label{caprace}
The following observation of P Caprace implies that the $\CAT (0)$ Conjecture holds for $v\geq 9$.  Let $(W,S)$ be a Coxeter system with $v(W)\geq 9$.   Then $W'$ is an irreducible Euclidean simplex reflection group of rank at least 10, and by inspection of the standard list of these groups, we find that $W$ decomposes as \hbox{$W=W' \times W_{S-S'}$}.  Thus $W$ acts  on $\EE^v$ via projection onto $W'$.  In particular, this action yields a decomposition of $W$ as a nonpositively curved $v$--simplex of special subgroups.  This decomposition is the $v$--splitting $\Lambda$ constructed above.
\end{remark}

In the next section, we will prove the $\CAT (0)$ Conjecture for $v=2$.  For $3\leq v\leq 8$, there are only finitely many different possibilities for $\Lambda '$ since, by the classification of irreducible Coxeter groups, there are only finitely many isomorphism types of subgroups $W'$.  These are the remaining open cases.

\subsection[Proof of the CAT(0) Conjecture in dimension 2]{Proof of the $\CAT (0)$ Conjecture in dimension 2} \label{Dim2Section}
We will now prove the $\CAT (0)$ Conjecture in the case $v=2$, ie, when the $v$--splitting of $W$ is a triangle of groups decomposition of $W$.  To do so, we will show that the Gersten--Stallings angles in this triangle of groups are the same as those for the natural decomposition of $W'$.  In particular, we will prove the following result (notation as in \fullref{GerStsection}).  Note that $\pi/\infty$ means 0.

\eject
\begin{prop} \label{Angles}
Suppose $(W,S)$ is a Coxeter system with distinct special subgroups $A$ and $B$.
\begin{enumerate}
\item \label{1} Let $m=\min\{ m_{ij}\st  s_i \in S_A-S_B,\ s_j \in S_B-S_A\}$.  Then,
\[\angle_W (A,B;A\cap B) =
\begin{cases}
0 & {\textrm{ if } A\subset B \textrm{ or } B \subset A};\\
\pi/{m} &\textrm{ otherwise.}\\
\end{cases}
\]
\item \label{2} If $C$ is a special subgroup of $W$ such that $A\cap C=\{1\}=B\cap C$, then \[\angle_{\langle A,B\rangle} (A,B;A\cap B ) = \angle_W (\langle A, C \rangle, \langle B, C \rangle ; \langle A\cap B ,C \rangle).\]
In particular, for $a\neq b \in S$ and $T \subset S-\{a,b\}$, we have \[\angle_W (\langle a,T\rangle, \langle b,T\rangle; \langle T \rangle) = \frac{\pi}{m_{ab}}.\]
\end{enumerate}
\end{prop}

Note that statement \eqref{2} of this proposition implies that the Gersten--Stallings angle at a vertex in a triangle of special subgroups of a Coxeter group does not change when (the same) additional Coxeter generators are added to every local group.

\begin{theorem} \label{dim2}
The \fullref{cat0conj} is true for $v=2$.
\end{theorem}

By \fullref{CAT0impliesFAn} above, this will complete the proof of \fullref{thm12}.

\begin{proof}[Proof of \fullref{dim2}]
As in \fullref{construct}, let $\Lambda$ be the triangle of groups decomposition of $W$ determined by $S' \subset S$, and let $X$ be the universal cover of $\Lambda$.  By statement \eqref{2} of \fullref{Angles}, the Gersten--Stallings angles of $\Lambda$ are the same as those of $\Lambda'$.  Since $W'$ is a hyperbolic or Euclidean triangle group, the sum of the angles of $\Lambda$ is therefore at most $\pi$.  Thus $X$ is $\CAT (0)$ by \fullref{CAT0triangles}.
\end{proof}

\begin{cor} \label{Pdim2}
The Coxeter $\FA _n$ Conjecture (\fullref{pdimconj}) is true for \hbox{$n=2$}.  In particular, a Coxeter group $W$ acts nontrivially on a $2$--dimensional complete connected $\CAT (0)$ space if and only if $W$ has an infinite special subgroup of rank $3$ with respect to some (hence any) Coxeter generating set.
\end{cor}

Note that this result gives a complete characterization of $\FA _n$ for many Coxeter groups, including all Coxeter groups of large, even and odd type.

If $\Lambda '$ is hyperbolic (ie, if $W'$ is a hyperbolic triangle group), then $X$ is piecewise-hyperbolic and hence $\CAT (-1)$.  Thus, we have also shown the following:

\begin{cor}
Suppose $(W,S)$ is a Coxeter system with $m_{ij}< \infty$ for all $i$ and $j$.  Suppose further that $W$ has a special subgroup of rank $3$ that is a hyperbolic triangle group.  Then $W$ acts nontrivially on a $\CAT(-1)$ $2$--complex.
\end{cor}

Before proceeding to the proof of \fullref{Angles}, we introduce the following:
\vspace{6pt}
\begin{notation} Let $(W,S)$ be a Coxeter system.
\begin{enumerate}
\item[(i)]
For $s,\ t \in S$ with $m_{st}< \infty$ we define $\alpha _{st}$ to be the alternating word in $s$ and $t$ of length $m_{st}$.  That is,
\begin{equation*}
\alpha_{st}=
\begin{cases}
\underbrace{st\ldots st}_{m_{st}} & \text{if $m_{st}$ is even};\\
\underbrace{st\ldots ts}_{m_{st}} & \text{if $m_{st}$ is odd}.
\end{cases}
\end{equation*}
\item[(ii)] We denote by $\widehat{w}$ the word $w$ with one letter removed.
\item[(iii)] To denote equality of words in the alphabet $S$, we will use the notation $=$ and for equality of the corresponding group elements in a group $H$ we will use $=_{H}$.
\end{enumerate}
\end{notation}

\vspace{6pt}
The following lemmas will be used in the proof of \fullref{Angles}.  The first is an immediate consequence of the Deletion Condition along with \fullref{uniquealtprod}.

\vspace{6pt}
\begin{lemma} \label{seesubword}
Let $(W,S)$ be a Coxeter system.  Suppose $s_1s_2\ldots s_k$ is a word such that for some $r,t \in S$ and $j\leq m_{rt}$, we have \[s_1s_2\ldots s_k=_W\underbrace{rtr\ldots}_{j}.\]  Then, for some indices $i_1<i_2< \ldots < i_j$, we have $s_{i_j}=r$ for $j$ odd and $s_{i_j}=t$ for $j$ even (or perhaps vice versa if $j=m_{rt}$).  That is, as a word in the alphabet $S$, we have $s_{i_1}s_{i_2}\ldots s_{i_{j}} = \underbrace{rtr\ldots}_j \text{ (or possibly $\displaystyle \underbrace{trt\ldots}_{j}$ if $j=m_{rt}$})$.
\end{lemma}
\vspace{6pt}

The notation used below is as in \fullref{XX'reduced}.
\vspace{6pt}

\begin{lemma} \label{findingXX'reduced}
Let $(W,S)$ be a Coxeter system.  Let $r,t, \in S$ and let $w(r,t)$ be an alternating word in $r$ and $t$ of length strictly less than $m_{rt}$.  Suppose $w$ is the element of $W$ represented by the word $w(r,t)$.  Let $i(w)$ (resp.\ $j(w)$) be the first (resp.\ last) letter in the word $w(r,t)$.  Then $w$ is \hbox{$(I-i(w), J-j(w))$--reduced} for all $I, J\subset S$ .
\end{lemma}
\vspace{6pt}

\begin{proof}[Proof of \fullref{findingXX'reduced}]
By \fullref{XX'reduced}, it suffices to show that for all $s \in
I-i(w)$ and all $s' \in J-j(w)$, we have $\ell (sw) > \ell(w)$ and
$\ell(ws') > \ell(w)$.  Note that the Deletion Condition implies that
$\ell(sw) \neq \ell(w)$.  Suppose that $\ell(sw) < \ell (w)$ for some
$s \in I -i(w)$.  By the Strong Exchange Condition
(\fullref{deletion}), we have $ w(r,t)=_Ws\widehat{w(r,t)}$.  However,
since $w(r,t)$ is the unique reduced word representing the element $w$
of $W$ (see \fullref{uniquealtprod}), we have
$w(r,t)=s\widehat{w(r,t)}$ as words in the alphabet $S$.  Thus $s=r$
or $s=t$.

Without loss of generality, we may assume $w(r,t)$ begins with $r$.  Now $s\neq r$ since $r=i(w)$ and $s\in I-i(w)$.  So $s=t$ and $sw(r,t)$ is an alternating product of $r$ and $t$ beginning with $t$ and of length $\ell(w) +1$.  Note that $sw(r,t)$ is reduced since $\ell(w) < m_{rt}$, so $\ell (sw)=\ell(w)+1 > \ell(w)$.  This is a contradiction.

The case of $s' \in J-j(w)$ is analogous.
\end{proof}

We will now prove \fullref{Angles}.

\begin{proof}[Proof of \fullref{Angles}]
First note that statement \eqref{2} follows immediately from statement \eqref{1} since  \fullref{specialchar} implies that $\langle A,C\rangle \cap \langle B, C \rangle = \langle A\cap B , C \rangle$ and $S_{\langle A, C\rangle} - S_{\langle B, C \rangle} =S_A-S_B$.  We proceed now to prove statement \eqref{1}.

Let $G = A \ast _{A\cap B} B $ and let $\rho \co G \twoheadrightarrow \langle A,B \rangle \subset W$ be the natural surjection.  If $A\cap B =A$ or $A\cap B=B$, then \hbox{$\angle_W (A,B;A\cap B) =0$} since the induced map $G\rightarrow \langle A, B \rangle$ is an isomorphism.  The same is true if $m_{ab} =\infty$ for all $a\in S_A-S_B$ and $b\in S_B-S_A$.   Assume otherwise.

By definition, $\angle_W(A,B;A\cap B) = \frac{\pi}{k}$ where $2k$ is the minimal normal form length among nontrivial elements in $\ker (\rho)$.  Let $a\in S_A-S_B$ and $b \in S_B-S_A$ such that $m_{ab}$ is finite.  Then, \hbox{$(ab)^{m_{ab}} =_W 1$} but $(ab)^{m_{ab}}$ is nontrivial in $G$ (see \fullref{amalgcox}).  Therefore $(ab)^{m_{ab}}$ is a nontrivial element in $\ker (\rho)$.  Moreover, by \fullref{nflred} (iii), the normal form length in $G$ of $(ab)^{m_{ab}}$ is $2m_{ab}$, so $k \leq m_{ab}$ for all such choices of $a$ and $b$.  To prove part (1), it thus suffices to show that there exist elements $a\in S_A-S_B$ and $b \in S_B -S_A$ such that $k\geq m_{ab}$.

Let $M$ be the Coxeter matrix of $(W,S)$ and, as described in
\fullref{amalgcox}, let $M'$ be the Coxeter matrix of $(G,S_A \cup
S_B)$.  Then every elementary $M'$--operation is an elementary
$M$--operation (see \fullref{elemdefn}).  Moreover, the only
elementary $M$--operations for elements in $\langle A, B \rangle$
which are not also elementary $M'$--operations are those which replace
$\alpha_{ab}$ by $\alpha_{ba}$ (or vice versa) for some $a \in S_A -
S_B$ and $ b \in S_B - S_A.$

Let $g\in \ker(\rho)$.  By \fullref{titswp}, since $\displaystyle g =_W 1$ there is an $M$--reduction of $g$ to the identity.  Since $g \neq _G 1$ we know that $g$ does not $M{'}$--reduce to the identity.  So, the $M$--reduction of $g$ to the identity in $W$ involves at least one elementary $M$--operation which is not also an elementary $M'$--operation.  In particular, for some $a \in S_A-S_B$ and $b \in S_B -S_A$ with $m_{ab}$ finite, some $M'$--reduced expression of $g$ must have a subword of the form $\alpha_{ab}$ or $\alpha_{ba}$.

Hence we may choose $g$ nontrivial in $ \ker (\rho)$ along with $a\in S_A - S_B$ and $b\in S_B-S_A$ so that $m_{ab}$ is minimal among all possible such choices also satisfying both of the following:
\begin{enumerate}
\item[(i)] The normal form length of $g$ in $G$ is $2k$.
\item[(ii)] Some $M'$--reduced expression of $g$ has $\alpha_{ab}$  or $\alpha_{ba}$ as a subword.
\end{enumerate}

Without loss of generality, we may assume that the alternating expression of (ii) begins with $a$, ie, that $\alpha_{ab}$ is a subword of an $M'$--reduced word representing $g$.  For notational convenience we will also use $g$ to denote this particular choice of $M'$--reduced word.

Let $m=m_{ab}$.  To complete the proof of part (1), we will show $k\geq m$.  We consider separately the cases when $m$ is even and odd.

\textbf{Case 1}\qua  Suppose $m$ is even.

Since $g$ is $M'$--reduced, by \fullref{coxred} (ii), we may assume
\[g=a_1b_2\ldots b_{j-1}a_j \alpha_{ab} b_{m+j-1}a_{m+j} \ldots a_{2k-1}b_{2k}\]
with $a_ja,a_i \in A-\{1\}$ for $ i \neq j$ and similarly $b b_{m+j-1},b_i \in B-\{1\}$ for $i \neq m+j-1$.

If $m=2$ then we need only show $k\neq 1$.  If $k=1$ then $g=a_1abb_2$ for some $a_1 \in A$ and $b_2 \in B$.  Then $a_1a=_Wb_2^{-1}b$ so $a_1a$ and $bb_2$ are elements of $A\cap B$ and thus $g\in A \cap B$.  In particular, the fact that $g=_W1$ implies that $g=_G1$ in $G$ since $\rho \mid _{A\cap B}$ is injective.  This is a contradiction as $g$ is a nontrivial element in $\ker (\rho) \subset G$.

Assume now that $m> 2$.  Conjugation in $G$ yields another nontrivial element $h$ of $\ker(\rho)$ represented by the word
\[
h= \alpha_{ab} b_{m+j-1} a_{m+j} \dots  a_{2k-1}b_{2k} a_1b_2 \ldots b_{j-1}a_j.\]

Since $h=_W1$ and $\alpha_{ba}=_W \alpha_{ab}^{-1}$, we have
\begin{equation} \label{eveninW}
\alpha_{ba} =_W b_{m+j-1} a_{m+j} \dots  a_{2k-1}b_{2k} a_1b_2 \ldots b_{j-1}a_j.
\end{equation}
As $a\notin S_B$, the generator $a$ cannot appear in any $M'$--reduced expression of $b_i$ for any $i$ by \fullref{specialchar} (and similarly for $b$ and $a_i$).  Applying \fullref{seesubword} we therefore find that $m \leq j + (2k-(m+j-1) +1) = 2k-m+2$ so $k\geq m-1$.  If $k>m-1$ then we are done, so we assume $k=m-1$.

Let $\tilde{h}$ denote the word on the right side of equation \eqref{eveninW}.  Since $k=m-1$, the word $\tilde{h}$ consists of an alternating sequence of exactly $m$ elements of $A$ and $B$, ie, $\alt (\tilde{h})=m$ in the notation of \fullref{amalgcox}.  Thus the normal form length of the group element $\tilde{h}$ in $G$ is at most $m$.

Since $\alpha_{ba}$ is $M$--reduced, \fullref{titswp} implies that there is an $M$--reduction of $\tilde{h}$ to $\alpha_{ba}$.  Since $h$ is nontrivial in $G$ and $h=\alpha_{ab}\tilde{h}$, it follows that $\tilde{h} \neq _G \alpha_{ba}$ as $\alpha_{ab}\alpha_{ba}=_G1$.  Thus, by the same argument as above, there is an $M'$--reduced form of $\tilde{h}$ containing a subword of the form $\alpha_{a'b'}$ or $\alpha_{b'a'}$ for some $a'\in S_A-S_B$ and $b'\in S_B-S_A$.  By \fullref{coxred} (iii), the normal form length of $\tilde{h}$ in $G$ is thus at least $m_{a'b'}$.  However $m_{a'b'}\geq m$ by the minimality assumption on $m$, so we find that $m=m_{a'b'}$.  Moreover, the normal form length of $\tilde{h}$ is exactly $m_{a'b'}=m$ and there exist elements $\tilde{a}\in A$ and $\tilde{b} \in B$ so that $\tilde{h}=_G\tilde{b} \alpha_{b'a'} \tilde{a}$.  It follows that equation \eqref{eveninW} can be rewritten as $\alpha_{ba}=_W\tilde{b} \alpha_{b'a'} \tilde{a}.$  By \fullref{seesubword}, since $m>2$ we have $a'=a$ and $b'=b$.  Thus, we find $\tilde{h}=_G \tilde{b}\alpha_{ba}\tilde{a}$.  Hence, equation \eqref{eveninW} can be rewritten as
 $\alpha_{ba}=_W \tilde{b} \alpha_{ba} \tilde{a}$.  It follows that
 \begin{equation}\label{eqn:final} b\tilde{b} b \underbrace{ab\cdots ab}_{m-2}=_W \underbrace{ab\cdots ab}_{m-2} a\tilde{a}^{-1} a.\end{equation}

    Since $b\tilde{b}b$ is a subword of $h$, it is $M'$--reduced so is also $M$--reduced.  Then \fullref{findingXX'reduced} implies that the left side of equation \eqref{eqn:final} is $M$--reduced.  So, by \fullref{titswp}, there is a sequence of elementary $M$--operations taking the right side to the left side.  By the nontriviality of $h$, equation \eqref{eqn:final} does not hold in $G$, so at least one operation must be used which is not an $M'$--operation.  Hence a finite number of $M'$--operations transforms the right side of equation \eqref{eqn:final} into a word $w'$ containing a subword of the form $\alpha_{a'b'}$ or $\alpha_{b'a'}$ for some $a'\in S_A-S_B$ and $b'\in S_B-S_A$.  Note that $m_{a'b'}\geq m$ by the minimality assumption on $m$.  Then, in the notation of \fullref{amalgcox}, we find \[\alt(w')\geq \alt(\alpha_{a'b'})=m_{a'b'}\geq m> \alt(\underbrace{ab\cdots ab}_{m-2} a\tilde{a}^{-1}a).\] This contradicts \fullref{coxred} (i) and thus completes the even case.

\textbf{Case 2}\qua  Suppose $m$ is odd.

By an argument analogous to the one in the even case, we find a nontrivial $M'$--reduced $h\in \ker (\rho)$ so that $h=_G\alpha_{ab} \tilde{a}_1 \alpha_{ab} \tilde{a}_2$
for some $\tilde{a}_1, \tilde{a}_2 \in A$.  This implies that \begin{equation}\label{eqn:odd} a\tilde{a_1} a \underbrace{ba\cdots ab}_{m-2}=_W \underbrace{ba\cdots ab}_{m-2} a\tilde{a_2}^{-1} a.\end{equation}

As in Case 1, the left side of equation \eqref{eqn:odd} is $M$--reduced by \fullref{findingXX'reduced} since $a\tilde{a_1} a$ is a subword of $h$.  Again there is a sequence of $M'$--operations that when applied to the right side of equation \eqref{eqn:odd} results in a subword $\alpha_{a'b'}$ or $\alpha_{b'a'}$ with $m_{a'b'}\geq m$.  This again contradicts \fullref{coxred} (i).
\end{proof}

\section[Maximal FA _n subgroups]{Maximal $\FA _n$ subgroups}\label{section:MaxFAn}
Let $W$ be a Coxeter group.  We say that a subgroup $H\subset W$ is {\em maximal $\FA _n$} if $H$ has property $\FA _n$ and $H$ is not properly contained in any other subgroup of $W$ with $\FA _n$.  Recall that every $(n+1)$--spherical special subgroup of $W$ (with respect to {\em any} Coxeter generating set of $W$) has $\FA _n$ by \fullref{Coxnfinite}.  We conjecture that all maximal $\FA _n$ subgroups of $W$ arise in this way.

\begin{FAnconj}
Let $(W,S)$ be a Coxeter system.  A subgroup $H \subset W$ is maximal $\FA _n$ if and only if $H = wAw^{-1}$ for some maximal $(n+1)$--spherical special subgroup $A$ of $W$ and $w \in W$.
\end{FAnconj}

Note that if the Coxeter $\FA _n$ Conjecture is true, then the above statement is equivalent to the formulation given in the Introduction.

The Maximal $\FA_1$ Conjecture was proven by Mihalik and Tschantz in \cite{MihalikTschantz}.  In this section, we use a modification of their arguments.  We prove the following reformulation of \fullref{maxFAnfromCAT0}.

\begin{theorem} \label{strongMaxFAn}
Suppose the $\CAT(0)$ Conjecture holds for all $v\leq n$.  Then the Maximal $\FA_n$ Conjecture holds for every Coxeter group.
\end{theorem}

By \fullref{thm12}, we will thus have the following:

\begin{cor}
The Maximal $\FA_2$ Conjecture is true for every Coxeter group.
\end{cor}

Because $\FA _n$ is a property of the group rather than of a particular presentation, we immediately conclude the following:

\begin{cor}
For all Coxeter groups $W$, and for all $n$ for which the Maximal $\FA_n$ Conjecture holds, the set of conjugates of maximal $(n+1)$--spherical special subgroups of $W$ is independent of the Coxeter presentation.  In particular, this is true for $n=1$ and $n=2$.
\end{cor}

The next proposition shows that the only candidates for maximal $\FA_n$ subgroups are indeed conjugates of $(n+1)$--spherical special subgroups.

\begin{prop} \label{FAncontainment}
Let $(W,S)$ be a Coxeter system.  Suppose a subgroup $H$ of $W$ has property $\FA_n$.  If the $\CAT (0)$ Conjecture holds for all $v\leq n$, then there is a $w\in W$ such that $H \subset wBw^{-1}$ for some $(n+1)$--spherical special subgroup $B$ of $W$.
\end{prop}

\begin{proof}
Let $v=\max\{m \st W \textrm{ is $m$--spherical}\}$ as in \fullref{construct}.  If $n < v$, then $W$ is $(n+1)$--spherical, and we are done.  So, we may assume $n \geq v$.  Note then that $H$ has $\FA_v$.  By our construction in \fullref{construct}, there is a special $v$--splitting $\Lambda$ of $W$ with vertex groups special subgroups of $W$ of rank $|S|-1$.  Since $H$ has $\FA_v$ and the $\CAT (0)$ Conjecture is assumed true for $v$, we know that $H$ fixes a point in the action of $W$ on the universal cover of $\Lambda$.  Therefore, $H \subset w_1 B_1 w_1^{-1}$ for some special subgroup $B_1 \subset W$ of rank $|S| - 1$ and some $w _1 \in W$.  If $B_1$ is $(n+1)$--spherical, then we are done.  Otherwise, let \hbox{$v_1=\max\{m \st B_1 \textrm{ is $m$--spherical}\}$.}  Then $v_1 \leq n$ and there is a special \hbox{$v_1$--splitting} of $B_1$.  As before, $H \subset w_2 B_2 w_2^{-1}$ where $B_2$ is a special subgroup of rank $|S| - 2$ and $w_2 \in W$.  Once again, if $B_2$ is $(n+1)$--spherical, then we are done.  Otherwise, we can continue splitting in this way.  Since $S$ is finite and the number of generators of $B_k$ is $|S|-k$, the process terminates.  So, we find $H \subset w_kB_kw_k^{-1}$ for some $(n+1)$--spherical $B_k$.
\end{proof}

\begin{remark}\label{nonmax} \fullref{FAncontainment} implies in particular that the only possible subgroups of $W$ which could be maximal $\FA_n$ are conjugates of $(n+1)$--spherical special subgroups.   Moreover, if $B \subset A$ are $(n+1)$--spherical subgroups of $W$, then $wB w^{-1} \subset wAw^{-1}$ for all $w \in W$, so if $H$ is a maximal $\FA _n$ subgroup of $W$, then $H$ is a conjugate of a {\em maximal} $(n+1)$--spherical subgroup.  Since all $(n+1)$--spherical subgroups have $\FA _n$, it remains only to show that all conjugates of maximal $(n+1)$--spherical subgroups are {\em maximal} $\FA_n$ subgroups.
\end{remark}

The lemma below follows from Proposition 5.5 of Deodhar \cite{Deodhar} together with \fullref{inclgen}.
\begin{lemma} \label{ABflagsub}
Let $(W,S)$ be a Coxeter system and suppose $A$ is a maximal $(n+1)$--spherical  special subgroup of $W$.  If $B$ is any $(n+1)$--spherical special subgroup of $W$ such that $A \subset wBw^{-1}$ for some $w\in W$, then $S_A = S_B$.
\end{lemma}

We now complete the proof of \fullref{strongMaxFAn}.

\begin{proof}[Proof of \fullref{strongMaxFAn}]
As described in \fullref{nonmax}, it remains only to prove that all conjugates of maximal $(n+1)$--spherical subgroups are maximal $\FA _n$ subgroups.
Suppose now that $A$ is a maximal $(n+1)$--spherical subgroup and $w\in W$.  To show that $A$ is maximal $\FA _n$, it suffices to show that $A$ is not properly contained in a conjugate of any $(n+1)$--spherical subgroup.

Suppose $A \subset wBw^{-1}$ for some $(n+1)$--spherical subgroup $B$.  We will show that $A=wBw^{-1}$.  By \fullref{ABflagsub}, we have $S_A = S_B$, so $A=B$ and $A \subset wAw^{-1}$.  By \fullref{inclgen}, this implies that $S_A = dS_A d^{-1}$ where $d$ is an $(S_A,S_A)$--reduced element of $AwA$.  Thus we find
\[ wAw^{-1} = dAd^{-1}
 =  d\langle S_A\rangle d^{-1}
 = \langle dS_Ad^{-1} \rangle
 =  \langle S_A \rangle
 =  A.\]
Hence $A=wAw^{-1}=wBw^{-1}$ so $A$ is indeed a maximal $\FA _n$ subgroup.  It immediately follows that the conjugates of $A$ are all maximal $\FA _n$ as well.
\end{proof}

\begin{remark}
Note that the above arguments also apply to maximal strong $\FA _n$ subgroups.  So, in particular, if the $\CAT (0)$ Conjecture holds, the maximal strong $\FA _n$ subgroups of $W$ are the same as its maximal $\FA _n$ subgroups.
\end{remark}

\section{Proper actions}
In this section, we briefly consider the case of proper actions of Coxeter groups.  The {\em $\CAT(0)$ dimension} of a group $G$, denoted {\em $\dim_{ss}(G)$}, is defined to be the minimal dimension of a complete $\CAT (0)$ space on which $G$ acts properly by semisimple isometries. $\CAT (0)$ dimension has been studied by Bridson \cite{Bridson01} and Brady and Crisp \cite{BradyCrisp02}, among others.

The following lemma is immediate from the definitions.
\begin{lemma}
Suppose $G$ is a group and $H$ is an infinite subgroup of $G$ with strong $\FA _n$.  Then $\dim _{ss} (G) > n$.
\end{lemma}

Together with \fullref{Coxnfinite}, this gives the following lower bound on the $\CAT (0)$ dimension of Coxeter groups.
\begin{cor} \label{properdim}
Let $(W,S)$ be a Coxeter system.  Then
\[
\dim _{ss}(W) > \max\{n \st W \textrm{ has an infinite } (n+1)\textrm{--spherical special subgroup} \}.\]
\end{cor}

In fact, we obtain a better lower bound by applying \fullref{leraygroupcriterion} to the set of finite special subgroups of $W$.  Let $L(W,S)$ denote the simplicial complex with vertices corresponding to the elements of $S$ and such that a subset $U\subset S$ spans a simplex in $L(W,S)$ if and only if $W_U$ is finite.  The complex $L(W,S)$ is called the {\em nerve} of the Coxeter system $(W,S)$.
\begin{cor}
Let $(W,S)$ be a Coxeter system.  Then
\[ \dim_{ss}(W) > \max\{k \st H_k(L(W_T,T))\neq 0 \mbox{ for some } T \subset S\}.\]

\end{cor}

\begin{proof}
(Notation as in \fullref{homimp}.) Let $T\subset S$.  Suppose $\phi\co W \rightarrow \Isom(X)$ gives an action of $W$ on a complete $\CAT (0)$ space $X$ via semisimple isometries.  Then the Bruhat--Tits Fixed Point Theorem (\fullref{BT}) implies that $L(W_T,T)\subset \cN (T, \phi)$.  On the other hand, \fullref{leraygroupcriterion} implies that $H_k(\cN(T,\phi))=0$ for $k\geq \dim(X)$.  So, if \hbox{$H_k(L(W_T,T))\neq 0$} for some $k\geq \dim(X)$, then \hbox{$L(W_T,T)\neq \cN(T,\phi)$}.  In particular, some infinite subgroup of $W$ must fix a point, so the action is not proper.
\end{proof}

On the other hand, a Coxeter group $W$ acts properly on its Davis--Moussong Complex $\Sigma_{DM}$ (see \fullref{Coxactions}).  It follows that $\dim (\Sigma _{DM})$ is an upper bound for $\dim_{ss}(W)$.  In particular,
\[\tag{$*$}
\dim_{ss} (W) \leq \max \{ n \st W\textrm{ has a finite special subgroup of rank }n\}. \label{davisbd}\]
A natural question then is how the bounds on $\dim _{ss}(W)$ given by the above inequalities are related.  The following example shows that the upper bound given by (\ref{davisbd}) is not optimal.

\begin{figure}[ht!]
\begin{minipage}{0.47\linewidth}
\hypertarget{double2}{\cl{\includegraphics{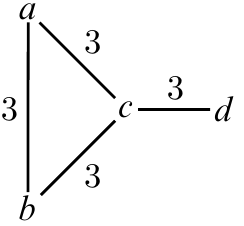}}}
\caption{\newline Coxeter diagram of $W$}\label{properexamplediagram}
\end{minipage}
\addtocounter{figure}{1}
\begin{minipage}{0.47\linewidth}
\cl{\includegraphics{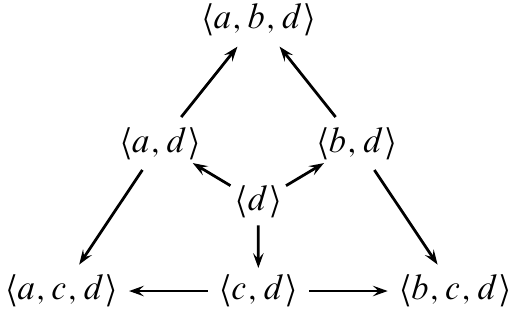}}
\caption{Triangle of\newline finite groups}\label{properexampletriangle}
\end{minipage}
\end{figure}
 
\begin{eg}
Consider the Coxeter group $W$ with Coxeter diagram given by \mbox{\fullref{properexamplediagram}}.  Note that $W$ is $2$--spherical but not $3$--spherical since $W_{\{a,b,c\}}$ is infinite.  So by \fullref{Pdim2} we know that $W$ has strong $\FA_1$ but not strong $\FA_2$.  In particular, the triangle of groups decomposition of $W$ in \hyperlink{double2}{Figure 7} gives a $\CAT (0)$ $2$--complex on which $W$ acts nontrivially by semisimple isometries.
Moreover, the point stabilizers under this action are conjugates of the local groups in \hyperlink{double2}{Figure 7}.  In particular, all point stabilizers are finite, so the action is proper.  Hence $\dim _{ss} (W) \leq 2$. Moreover, by \fullref{properdim}, since $W_{\{a,b,c\}}$ is infinite and $2$--spherical, we have $\dim_{ss} (W) \geq 2$.  So, the $\CAT (0)$ dimension of $W$ is $2$.  However, the Davis--Moussong complex $\Sigma _{DM}$ is $3$--dimensional since, for example, $W_{\{a,b,d\}}$ is finite.
\end{eg}

\bibliographystyle{gtart}
\bibliography{link}

\end{document}